\ifx\LocalElevenPtArticleFormatLoaded\undefined
\documentclass[11pt]{article}
\usepackage{theorem} 
\usepackage{upref,amsmath,amscd,amssymb} 
\usepackage{latexsym} 
\usepackage{diagrams}
\diagramstyle[small,labelstyle=\scriptstyle,midshaft,PS]
\usepackage[german,english]{babel} 
\usepackage{epsfig}
\usepackage{xypic}
\fi
\newif\ifpdf\ifx\pdfoutput\undefined\pdffalse\else\pdfoutput=1
\usepackage[pdftitle={A note on the Intersection of two Veronese surfaces}, 
         colorlinks=true,
         hypertexnames=false,
         linkcolor=blue,
         citecolor=red,
         pdfstartview=XYZ,
         pdfpagemode=FitWidth,pdftex]{hyperref}
\usepackage[nottoc]{tocbibind}
\usepackage[pdftex]{thumbpdf} 
\usepackage[pdftex]{color}\pdftrue\fi

\newtheorem{theorem}{Theorem}[section]
\usepackage{theorem}

\newtheorem{lemma}[theorem]{Lemma}

\newtheorem{proposition}[theorem]{Proposition}
\theorembodyfont{\rm}
\newtheorem{remark}{Remark}

\newcommand{\PP}{{\mathbb{P}}}
\newcommand{\QQ}{{\mathbb{Q}}}

\newcommand{\ZZ}{{\mathbb{Z}}}

\newcommand{\calO}{{\cal O}}

\newcommand{\calI}{{\cal I}}

\newcommand{\calE}{{\cal E}}

\newcommand{\Sec}{{\operatorname{Sec}}}
\newcommand{\Sym}{{\operatorname{Sym}}}

\newenvironment{Proof}{\begin{ProofwCaption}{Proof}}{\end{ProofwCaption}}
\newenvironment{Proof*}[1]{\begin{ProofwCaption}{{#1}}}{\end{ProofwCaption}}
\newenvironment{ProofwCaption}[1]%
  {\addvspace\theorempreskipamount \noindent{\it #1.}\rm}%
  {\qed \par \addvspace\theorempostskipamount}
\newcommand{\qedsymbol}{\mbox{$\Box$}}
\newcommand{\qed}{\quad\qedsymbol}
\begin{document}

\title{A note on the Intersection of Veronese Surfaces}
\author{D.~Eisenbud, K.~Hulek and S.~Popescu }
\date{February 10, 2003}
\maketitle

\section{Introduction}
The main purpose of this note it to prove the following

\begin{theorem}\label{theo0.1}
Any two Veronese surfaces in $\PP^5$ whose intersection is zero-dimensional
meet in at most $10$ points (counted with multiplicity).
\end{theorem}

Our initial motivation for this note comes from our paper \cite{EGHPO}
where we study linear syzygies of homogeneous ideals generated by
quadrics and their restriction to subvarieties of the ambient
projective space with known (linear) minimal free resolution.  A
direct application of the techniques in \cite[Section 3]{EGHPO} 
shows that the homogeneous ideal of a zero-dimensional intersection of two
Veronese surfaces in $\PP^5$ is $5$-regular (see also Lemma \ref{five-reg}
below), which  yields only an upper bound of $12$ for its degree,
cf. Section \ref{Reduction-step}.

Section \ref{Veronese-on-hyperquadrics} analyzes Veronese surfaces on
hyperquadrics. The observation that two Veronese surfaces on a smooth hyperquadric
$Q\subset\PP^5$ meeting in a zero-dimensional subscheme, must meet in a
subscheme of length $10$ or $6$ is classical and goes back to Kummer
\cite{Kummer} and Reye \cite{Rey} (see also \cite{Jes1} for historical
comments): By regarding the smooth hyperquadric $Q\subset\PP^5$ as the
Pl\"ucker embedding of the Grassmannian of lines
${\textrm{Gr}}(\PP^1,\PP^3)$, a Veronese surface on $Q\subset\PP^5$
is, up to duality, the congruence of secant lines to a twisted cubic
curve and thus has bidegree $(1,3)$. More precisely, the congruence
has one line passing through a generic point of $\PP^3$ and $3$ lines
contained in a generic plane. Thus Schubert calculus yields that the
possible intersection numbers of two Veronese surfaces on $Q$ are
either $10$ or $6$. See Proposition \ref{prop2.1} below and the
following remark, or the computation of the number of common chords of
two space curves in \cite[page 297]{GH}.

The case of two Veronese surfaces in $\PP^5$ meeting in $10$ simple
points has also been investigated in relation with association
(projective Gale transform) by Coble \cite{Coble1}, Conner \cite{Con}
and others. In \cite[Theorem 26]{Coble1} Coble claims that $10$ points
in $\PP^5$ which are associated to the $10$ nodes of a symmetroid in
$\PP^3$, the quartic surface defined by the determinant of a symmetric
$4\times 4$ matrix with linear entries in $\PP^3$, are the (simple)
intersection points of two Veronese surfaces in $\PP^5$  (see
Proposition \ref{Coble} below).  This is based on Reye's observation
\cite[page 78-79]{Rey} that $4\times 4$ symmetric matrices with linear
entries in $\PP^3$ are actually catalecticant with respect to suitable
bases and on the analysis in \cite {EiPo} of the Gale Transform of
zero-dimensional determinantal schemes. Section \ref{catalecticants}
contains a modern account of these results.

In Section \ref{How-many} we briefly discuss which intersection numbers 
$\le 10$ can actually occur for two Veronese surfaces in $\PP^5$ 
and in which geometric situation this can happen. For instance, we show that 
two Veronese surfaces in $\PP^5$ cannot intersect transversally in
9 points, however they may intersect in non-reduced zero-dimensional
schemes of this degree.

With the exception of Section \ref{catalecticants} all other results in 
this note are valid in arbitrary characteristic.

\section {A reduction step}
\label{Reduction-step}
We shall make essential use of the following lemma whose proof is
reminiscent of the linear syzygies techniques used in 
\cite[Section 3]{EGHPO}.

\begin{lemma}\label{five-reg}
If $X_1$ and $X_2$ are two Veronese surfaces in $\PP^5$ meeting in a
zero-dimensional scheme $W$, then the ideal sheaf 
${\cal I}_{W,\PP^2}$ of $W$ regarded as a subscheme of $\PP^2$ is
$5$-regular.
\end{lemma}

\begin{Proof}
The claim is equivalent to the vanishing $h^1({\cal I}_{W,\PP^2}(4))=0$. 
In order to see this we consider the minimal resolution of ${\cal I}_{X_1}$ in
$\PP^5$ and restrict it to $X_2$. This yields  a complex abutting to
${\cal I}_{W}$. Since $X_1$ has property $N_p$ for all $p$ one immediately
computes that $h^1({\cal I}_{W, \PP^2}(4))=0$. \hfill
\end{Proof}

As immediate consequence we get a first bound for the number of points where 
two Veronese surfaces whose intersection is zero-dimensional can meet.

\begin{proposition}\label{prop1.3}
Two Veronese surfaces in $\PP^5$ whose intersection is zero-dimensional meet
in at most $12$ points.
\end{proposition}
\begin{Proof}
Let $X_1$ and $X_2$ be two Veronese surfaces meeting in a
zero-dimensional scheme $W$, and let $d=\operatorname{length}(W)$.  By
Lemma \ref{five-reg}, $W$ regarded as a subscheme of $X_2\cong \PP^2$
imposes independent conditions on plane quartics, in particular if
$d\le 15$. On the other hand $X_2\subset\PP^5$ is cut out by quadrics
scheme-theoretically, and thus the quartics in $H^0({\cal
I}_{W,\PP^2}(4))$ must cut out $W\subset\PP^2$ scheme-theoretically
too -- in particular there are at least two. If there were only two,
then they would form a complete intersection, generating a saturated
ideal, and thus $W\subset\PP^2$ would be a complete intersection of
two plane quartics, which cannot be $5$-regular.  It follows that
$h^0({\cal I}_{W,\PP^2}(4))\ge 3$, and thus that $d\le 12$.\hfill
\end{Proof}

Actually the above proof yields also the following estimate

\begin{proposition}\label{prop1.4}
Let $X_1$ and $X_2$ be two Veronese surfaces in $\PP^5$ meeting in a
zero-dimensional scheme $W$ of length $d$ with $10 \le d \le 12$.
Then $X_1\cup X_2$ lies on at least $d-9$ quadrics.
\end{proposition}

\begin{Proof}
Let $a=h^0({\cal I}_{X_1\cup X_2}(2))$. Then on one hand
$h^0({\cal I}_{W,\PP^2}(4)) \ge h^0({\cal I}_{X_2}(2))-a=6-a$,
on the other hand, by Lemma \ref{five-reg}, we know that 
$h^0({\cal I}_{W,\PP^2}(4))=15-d$. Combining the two 
proves the claim of the proposition.
\hfill\end{Proof}

\section{Veronese surfaces on hyperquadrics}
\label{Veronese-on-hyperquadrics}
In this section we analyze the intersection of two Veronese surfaces
which meet in finitely many points, in the case where the two surfaces
lie on a common hyperquadric, resp. a pencil of hyperquadrics. We begin
with the classical case of congruences of lines:

\begin{proposition}\label{prop2.1}
Assume that $X_1$ and $X_2$ are Veronese surfaces which meet in finitely many
points and assume moreover that there exists a smooth quadric hypersurface $Q$
with $X_1\cup X_2\subset Q$. Then $X_1.X_2=6$ or $10$.
\end{proposition}

\begin{Proof}
Since $Q$ is smooth it is isomorphic to the Grassmannian
${\textrm{Gr}}(\PP^1,\PP^3)$ and it is well known that
$H^4({\textrm{Gr}}(\PP^1,\PP^3),\ZZ)=\ZZ\alpha+\ZZ\beta$ where
$\alpha$ and $\beta$ are $2$-planes. It follows from the double point
formula (see for instance \cite{HS}) that every Veronese surface on
${\textrm{Gr}}(\PP^1,\PP^3)$ has class $3\alpha+\beta$ or
$\alpha+3\beta$. Since $\alpha^2=\beta^2=1$ and $\alpha\beta=0$ it
follows that $X_1X_2=10$ or $6$ depending on whether $X_1$ and $X_2$
belong to the same class or not.  \hfill \end{Proof}

\begin{remark}
As already mentioned in the introduction, a Veronese surface of class
$3\alpha+\beta$ on the Grassmannian ${\textrm{Gr}}(\PP^1,\PP^3)\cong
Q\subset\PP^5$ is the congruence of secant lines to a twisted cubic
curve (where $\beta$ is the cycle of lines passing through a point of
$\PP^3$ while $\alpha$ is the cycle of lines in a plane). Passing to
the dual $\PP^3$ exchanges $\alpha$ and $\beta$, so up to duality the
same construction accounts also for Veronese surfaces of class
$\alpha+3\beta$.  It is easy to see, for instance by using Kleiman's
transversality theorem, that both cases described in Proposition
\ref{prop2.1} actually occur.
\end{remark}

\begin{proposition}\label{prop2.2}
Assume that $X_1$ and $X_2$ are Veronese surfaces in $\PP^5$ which
meet in finitely many points and assume that there exists a rank $5$
hyperquadric $Q$ containing both $X_1$ and $X_2$. Then $X_1.X_2=8$.
\end{proposition}

\begin{Proof}
We first claim that $X_1$ and $X_2$ do not pass through the vertex $P$ of 
the quadric cone $Q$. Otherwise projection from $P$ would map the Veronese surface to a cubic scroll
contained in a smooth hyperquadric $Q'\subset\PP^4$. But this is impossible, since
by the Lefschetz theorem every surface on $Q'$ is a complete intersection and
hence has even degree.
Blowing up the point $P$ we obtain a diagram
$$
\diagram
\tilde Q\rto^p\dto^\pi &Q\\
Q'
\enddiagram
$$ 
where $\pi$ gives $\tilde Q$ the structure of a $\PP^1$-bundle over
$Q'$.  Since $X_1$ and $X_2$ do not go through the point $P$ they are
not blown up and we will, by abuse of notation, also denote their
pre-images in $\tilde Q$ by $X_1$ and $X_2$.  The Chow ring of $\tilde
Q$ is generated by $H=p^{\ast}(H_{\PP^5})$ and
$H'=\pi^{\ast}(H_{\PP^4})$. Clearly $H^4=H^3 H'=H^2(H')^2=H(H')^3=2$
and $(H')^4=0$. Let $E$ be the exceptional divisor of the map
$p$. Then $E=\alpha H+\beta H'$ and from $EH^3=0$ and $E(H')^3=2$ one
deduces $\alpha=1$ and $\beta=-1$, i.e. $E=H-H'$. The surfaces $X_i$
have class
$$
X_i=\alpha_i H^2+\beta_i H H'+\gamma_i(H')^2.
$$
{}From $X_i H^2=\deg X_i=4$ one computes $\alpha_i+\beta_i+\gamma_i=2$. Since
$X_i$ does not meet $E$ we have $X_iEH'=0$ and from this one deduces
$\gamma_i=0$ and hence
$$
X_i=\alpha_i H^2+(2-\alpha_i) H H'.
$$
But then $X^2_i=8$ and this proves the claim.
\hfill\end{Proof}

We analyze next what happens when $X_1$ and $X_2$ lie on a pencil of
hyperquadrics.

\begin{proposition}\label{prop2.3}
Assume that the Veronese surfaces $X_1$ and $X_2$ meet in finitely many
points and assume that they are contained in a pencil of hyperquadrics
$\{\lambda_1 Q_1+\lambda_2 Q_2=0\}$. Then for a general hyperquadric $Q$ in
this pencil
$$
X_1\cap X_2\cap\operatorname{Sing } Q=\emptyset.
$$
\end{proposition}

\begin{Proof}
Assume that this is false. Since $X_1\cap X_2$ is a finite set it
follows that there exists some point $P\in X_1\cap X_2$ which is
singular for every hyperquadric $Q$ in this pencil. {}From what we saw in
the proof of Proposition \ref{prop2.2} this also shows that the
general quadric in this pencil has rank at most $4$ (and at least $3$
since both surfaces $X_i$ are non-degenerate). Projecting from $P$
maps $X_1$ and $X_2$ to rational cubic scrolls $Y_1$ and $Y_2$ in $\PP^4$,
respectively.  These cubic scrolls are contained in a pencil of
non-degenerate hyperquadrics $\{\lambda_1 Q'_1+\lambda_2 Q'_2=0\}$ whose general
member has rank $3$ or $4$. For degree reasons this implies
$Y_1=Y_2$. (Incidentally this also shows that $X_1$ and $X_2$ are
contained in a net of hyperquadrics whose general element has rank $4$.)

Let $Y$ be the cone over $Y_1=Y_2$ with vertex $P$. We obviously have $X_1,
X_2\subset Y$. We now blow up in $P$ and obtain a diagram
$$
\diagram
\tilde Y\rto^{p\quad }\dto^\pi &Y\subset\PP^5\\
Y_1
\enddiagram
$$
where $\tilde Y$ is a $\PP^1$-bundle over $Y_1$. The Picard group of $Y_1$ is
generated by two elements $C_0$ and $F$ with $C^2_0=-1, C_0 F=1$ and $F^2=0$.
Let $F_1=\pi^{\ast}C_0$ and $F_2=\pi^{\ast}F_1$. Then the Chow group on $\tilde
Y$ is generated by $H=p^{\ast}(H_{\PP^5}), F_1$ and $F_2$.

For geometric reasons $H^3=3, F^2_1 H=-1, F_1 F_2 H=F_1 H^2=F_2 H^2=1$
and $F^2_2 H=F^2_1 F_2=F^2_2 F_1=0$. Let $E$ be the exceptional locus
of $p$ and let $\tilde X_i$ denote the strict transforms of the
Veronese surfaces. Then $\pi$ restricted to $\tilde X_i$ defines
isomorphisms between $\tilde X_i$ and $Y_1$. Since $X_1$ and $X_2$
intersect in only finitely many points and since both are contained in
the 3-dimensional cone $Y$ it follows that $X_1\cap X_2=\{P\}$, and
from this one concludes that $\tilde X_1\cap\tilde X_2=L$ where
$L=E\cap F_1$ is a projective line, respectively $\tilde X_1 \tilde X_2=aL$
for some $a\ge 1$. Next we want to determine the class of $\tilde X_i$
in $\tilde Y$.  Since the $\tilde X_i$ are sections of the
$\PP^1$-bundle $\pi:\tilde Y\rightarrow Y$ we find $\tilde
X_i=H+\beta_i F_1+\gamma_i F_2; i= 1,2$. Restricting this to $E$ and
using that $H$ is trivial on $E$ we immediately find that $\beta_i=1$
and $\gamma_i=0$, i.e. $\tilde X_i=H+F_1$. But then $\tilde X_1 \tilde
X_2=H^2+2 H F_1+F^2_1\neq aL$ where the latter inequality can be seen
e.g. by intersecting with $H$. This is a contradiction and the
proposition is proved.  \hfill\end{Proof}

\begin{proposition}\label{prop2.4}
Let $X_1, X_2$ be two Veronese surfaces in $\PP^5$ intersecting in a
finite number of points. If $X_1$ and $X_2$ are contained in a pencil
of hyperquadrics, then $X_1.X_2\le 10$.
\end{proposition}

\begin{Proof}
Let $r$ be the rank of a general element of this pencil of hyperquadrics. If $r=6$
or $5$ then the assertion follows from Proposition \ref{prop2.1}, or from
Proposition \ref{prop2.2}, respectively. On the other hand, since 
the surfaces $X_i$ are non-degenerate we must have $r\ge 3$. 

We shall first treat the case $r=4$. According to 
Proposition \ref{prop2.3} we can then choose a rank $4$ hyperquadric $Q$ with
$X_1\cap X_2\cap\operatorname{Sing} Q=\emptyset$.
Blowing up the singular line $L$ of $Q$ we obtain a diagram
$$
\diagram
\tilde Q\rto^{p}\dto^\pi &Q\subset\PP^5\\
\quad\qquad\qquad Q':=\PP^1\times\PP^1
\enddiagram
$$
where $\pi$ is the structure map of a $\PP^2$-bundle. We denote the
strict transforms of $X_i$ by $\tilde X_i$, $i=1,2$. Since $X_1\cap
X_2\cap\operatorname{Sing} Q=\emptyset$ we have $\tilde X_1.\tilde
X_2=X_1.X_2$.  Let $H=p^{\ast}(H_{\PP^5})$, let $L_1, L_2$ denote
the rulings of $\PP^1\times \PP^1$ and set $F_i=\pi^{\ast}L_i$,
$i=1,2$. Then $H^4=2, F_1 H^3=F_2 H^3=F_1 F_2 H^2=1$ and
$F^2_1=F^2_2=0$. Let $E$ be the exceptional locus of $p$. Its class
must be of the form $E=H+\alpha_1F_1+\alpha_2 F_2$. {}From $EF_1
H^2=H^2=EF_2H^2=0$ one deduces $\alpha_1=\alpha_2=-1$, i.e.
$E=H-F_1-F_2$.

Now let $X$ be any Veronese surface on $Q$. We want to determine the
possible classes of the strict transform $\tilde X$ of $X$ in $\tilde
Q$. Let $\tilde X=\alpha H^2+\beta_1F_1H+\beta_2F_2H+\gamma
F_1F_2$. {}From $\tilde X H^2=4$ we obtain
$2\alpha+\beta_1+\beta_2+\gamma=4$. A priory the singular line $L$ can
either be disjoint from $X$, meet it transversally in one point, be a
proper secant or a tangent of $X$. Projection from $L$ shows that only
the first and the third of these possibilities can occur. 

Assume first that $L$ and $X$ are disjoint. Then $\pi_{\mid\tilde X}:\tilde
X\rightarrow\PP^1\times\PP^1$ is a $2:1$ map which shows
$\alpha=2$. {}From $\tilde X EF_1=\tilde XEF_2=0$ we conclude
$\beta_1=\beta_2=0$ and hence $\gamma=0$, i.e. $\tilde X=2H^2$. 

Assume now that $L$ is a proper secant of $X$. Blowing up $Q$ along
$L$ then blows up $X$ in 2 points and the corresponding exceptional
curves are mapped to different rulings in $\PP^1\times \PP^1$. The map
$\pi_{\mid\tilde X}:\tilde X\rightarrow\PP^1\times\PP^1$ is now
birational and hence $\alpha=1$. From what we have just said it
follows that ${\tilde X}EF_1={\tilde X}EF_2=1$ and hence this implies
$\beta_1=\beta_2=1$. But then $\gamma=0$ and $\tilde X=H^2+HF_1+HF_2$.

Let $c_1=2H^2$ and $c_2=H^2+HF_1+HF_2$. The claim of the proposition
now follows for the rank $4$ case since $c^2_1=c^2_2=c_1c_2=8$.

It remains to deal with the case in which the general hyperquadric $Q$ in
the pencil has rank $3$. By Proposition \ref{prop2.3} we can again
assume that $X_1\cap X_2\cap\operatorname{Sing} Q=\emptyset$.
Projection from the singular locus of $Q$ gives a diagram
$$
\diagram
\tilde Q\rto^{p}\dto^\pi &Q\subset\PP^5\\
\qquad C\cong\PP^1
\enddiagram
$$
where $C$ is a conic section and $\pi:\tilde Q\rightarrow C$ is a
$\PP^3$-bundle. We denote $H=p^{\ast}(H_{\PP^5})$ and $F=\pi^{\ast}(pt)$. Then
$H^4=2, H^3 F=1$ and $F^2=0$. Let $X$ be any Veronese surface on $Q$ and denote
its strict transform on $\tilde Q$ by $\tilde X$. We first note that
$X\cap\operatorname{Sing} Q$ is a finite set. Otherwise
$X\cap\operatorname{Sing} Q$ would have to be a conic section and projection
from $\operatorname{Sing} Q$ would map $X$ onto a plane, not to a conic.
Finally let $E$ be the exceptional locus of $p$. The class of $E$ must be of
the form $E=H+\gamma F$ and from $E H^3=0$ it follows that $\gamma=-1$, i.e.
$E=H-F$. Now put $\tilde X=\alpha H^2+\beta H.F$. {}From $\tilde X H^2=4$ one
deduces that $2\alpha+\beta=4$. Since $X\cap\operatorname{Sing} Q$ is finite
one must have that $\tilde X EH=0$ and hence $\beta=0$. This shows that $\tilde
X=2H^2$ and the claim of the proposition follows since 
$X_1.X_2=\tilde X_1.\tilde X_2=4 H^4=8$.
\hfill\end{Proof}

\begin{remark}
The above proof shows that if the general element in the pencil of hyperquadrics containing
$X_1$ and $X_2$ has rank $3$ or $4$, then $X_1.X_2=8$.
\end{remark}

\section{Catalecticant symmetroids and Veronese surfaces}
\label{catalecticants}
In this section we prove Coble's claim \cite[Theorem 26]{Coble1}, 
mentioned in the introduction, that ten points in $\PP^5$ which are 
the Gale transform of the  nodes of a general  quartic 
symmetroid in $\PP^3$ are  the simple intersection points of two Veronese surfaces.

A {\sl quartic symmetroid} is the quartic surface in $\PP^3$ defined by the determinant 
of a  symmetric $4\times 4$ matrix with linear entries in $\PP^3$; for
general choices (of the matrix) the symmetroid has only ordinary double points as 
singularities (nodes) and their number is $10$, by Porteous' formula. 
These surfaces are sometimes called {\sl Cayley symmetroids}, as
Cayley initiated their study in \cite{Cay} (cf. \cite{Jes2}, but 
see \cite{Cos} for a modern account of Cayley's results and much more). 

A symmetric matrix whose diagonals are constant is called a
{\sl catalecticant} matrix.  Surprisingly enough, it turns out that a
symmetric $4\times 4$ matrix with linear entries in $\PP^3$ can always
be reduced to a catalecticant form (with respect to suitable
bases). This fact goes back to Reye \cite[page 78-79]{Rey} and Conner
\cite[page 39]{Con} and is (re)-proved below.

We will make use of the perfect pairing, called {\sl apolarity},
between forms of degree $n$ and homogeneous differential operators of
order $n$ induced by the action of $T=k[\partial_0,\ldots,\partial_r]$
on $S=k[x_0,\ldots,x_r]$ via differentiation:
$$\partial^{\alpha}(x^{\beta}) =  \alpha!\binom{\beta}{\alpha}x^{\beta-\alpha},$$
if $\beta \geq \alpha$ and $0$ otherwise, and where $\alpha$ and $\beta$ are multi-indices,
$\binom{\beta}{\alpha}=\prod \binom{\beta_i}{\alpha_i}$, and $k$ is a field of 
characteristic zero.

\begin{proposition} 
The Hessian matrix of a web of quadrics in $\PP^3$ is catalecticant
(with respect to a suitable basis) if and only if the quadrics in the
web annihilate the quadrics of a twisted cubic curve. (One says in
this situation that the web is ``orthic'' to the twisted cubic
curve.)
\end{proposition}

\begin{Proof} Let $q:W^*\rInto \Sym_2 V$ be the web of quadrics on $\PP^3=\PP(V)$. 
A twisted cubic $C\subset\check\PP^3=\PP(V^*)$ is defined by its
quadrics $H^0(\check\PP^3,\calI_C(2))$.  In suitable coordinates, say
$\partial_0,\ldots,\partial_3$, these are the minors of the matrix
$$
\begin{pmatrix}
\partial_0 & \partial_1 & \partial_2\\
\partial_1 & \partial_2 & \partial_3\\
\end{pmatrix}.
$$ 
In terms of the dual coordinates, $x_0,\ldots,x_3$ of $\PP(V)$, the web $q$ 
has the form 
\begin{align*}
a_0x_0^2+a_4x_1^2+a_7x_2^2+a_9x_3^2&+
2a_1x_0x_1+2a_2x_0x_2+2a_3x_0x_3\\
&+2a_5x_1x_2+2a_6x_1x_3+2a_8x_2x_3
\end{align*}
where $a_0,a_1,\ldots,a_9$ are linear forms in the variables of $W$.
Direct computation shows that a quadric in the web $q$ is annihilated
by the equations $\partial_0\partial_2-\partial_1^2,
\partial_0\partial_3-\partial_1\partial_2,
\partial_1\partial_3-\partial_2^2$ if and only if $a_2=a_4, a_3=a_5,
a_6=a_7$.  It follows that the web $q$ is orthic to the twisted cubic
$C$ iff its Hessian matrix has shape
$$
\begin{pmatrix}
b_0 & b_1 & b_2 & b_3\\
b_1 & b_2 & b_3 & b_4\\
b_2 & b_3 & b_4 & b_5\\
b_3 & b_4 & b_5 & b_6\\
\end{pmatrix}
$$
where $b_0,b_1,\ldots,b_6$ are linear forms in the variables of $W$, i.e.
it is catalecticant.\hfill\end{Proof}

It actually turns out that a $4\times 4$ symmetric matrix with linear
entries in $\PP^3$ can be represented in two different ways as a catalecticant 
matrix. Namely

\begin{proposition}\label{Reye}
There are exactly two twisted cubic curves whose defining quadrics are
annihilated by a general web of quadrics in $\PP^3$.
\end{proposition}

\begin{Proof} We use the same notation as in the proof of the previous proposition.
Namely, let $q:W^*\rInto \Sym_2 V$ denote a general web of quadrics in
$\PP^3=\PP(V)$, and let $\check\PP^3=\PP(V^*)$ denote the dual
space. We also choose coordinates as above such that
$S=k[x_0,\ldots,x_r]$ and $T=k[\partial_0,\ldots,\partial_r]$ are the
coordinate rings of $\PP(V)$ and $\PP(V^*)$, respectively. Let now
$q':U\subset W^*\rInto \Sym_2(V)$ be a general subnet.  The variety
$H(q')$ of twisted cubics in $\PP^(V^*)$ whose defining quadratic
equations are annihilated by the net $q'$ is the geometric realization
of a prime Fano threefold $X_q$ of genus $12$ (see \cite{Muk1} and
\cite{Sch}). Via the apolarity pairing, $X_q\subset
{\textrm{Gr}}(\PP^2,\PP(U_q))$, where the annihilator $U_q={(q^\perp)}_2\subset\Sym_2
V^*$ is a 7-dimensional vector space. As a subvariety of the
Grassmannian the Fano threefold $X_q$ is the (codimension 9) common
zero-locus of three sections of $\wedge^2 \calE$, where $\calE$ is the
dual of the tautological subbundle on ${\textrm{Gr}}(\PP^2,\PP(U_q))$,
cf. \cite{Muk1} or see \cite[Theorem 5.1]{Sch} for a complete proof.

The choice of a (general) subnet $q'$ is equivalent to the choice of a
(general) global section of $\calE$. But $\calE$ is a globally
generated rank $3$ vector bundle whose restriction $E=\calE_{\mid
X_q}$ has third Chern number $2$, as an easy direct computation
shows. By Kleiman's transversality theorem the (general) section of
$E$ corresponding to $q'$ must vanish exactly at two (simple) points
of $X_q$. These in turn correspond to two twisted cubic curves each of
whose defining quadrics are annihilated not only by the net $q'$, but
by the whole web $q$ (the zero locus of a section in $\calE$ is the
special Schubert cycle of subspaces lying in the hyperplane dual to
the section). This concludes the proof.\hfill\end{Proof}

Let $C\subset\PP^6$ be a rational normal sextic curve, and let
$S=\Sec(C)\subset\PP^6$ be its secant variety.  $S$ has degree 10,
since this is the number of nodes of a general projection of $C$ to a
plane.  The homogeneous ideal of $C$ is generated by the $2\times
2$-minors of either a $3\times 5$ or a $4\times 4$ catalecticant
matrix with linear entries, induced by splittings of
$\calO_{\PP^1}(6)$ as a tensor product of two line bundles of strictly
positive degree. Furthermore, it is known that the homogeneous ideal
of $S=\Sec(C)$ is generated by the $3\times 3$ minors of either of the
above two catalecticant matrices (see \cite{GP} or \cite{EKS}).

For $\Pi=\PP^3\subset\PP^6$ a general $3$-dimensional linear subspace,
the linear section $\Gamma=\Sec(C)\cap\Pi$ consists of $10$ simple
points in $\PP^3$ defined by the $3\times 3$ minors of a $4\times 4$
symmetric (even catalecticant) matrix with linear entries in the
variables of $\Pi$. Conversely, by Proposition \ref{Reye} above, the
set $\Gamma\subset\PP^3$ of $10$ nodes of a general quartic symmetroid in 
$\PP^3$ arises always as a linear section of the secant variety of the rational 
normal curve in $\PP^6$ (in two different ways).  Moreover, it follows that
$\Gamma$  can also be defined by the $3\times 3$-minors of each 
of two different $3\times 5$ catalecticant matrices with linear entries in $\PP^3$. 
Since the $2\times 2$ minors of these catalecticant matrices generate an
irrelevant ideal, we may apply \cite[Theorem 6.1]{EiPo} (see also
\cite[Example 6.3]{EiPo} for more details) to obtain the following

\begin{proposition} (Coble)
\label{Coble}
The Gale transform of the $10$ nodes of a general quartic symmetroid
in $\PP^3$ are the points of intersection of two Veronese surfaces in
$\PP^5$.
\end{proposition}

\begin{remark} 
1) A more careful analysis of the preceding argument shows that 
the needed generality assumptions on the quartic symmetroid are satisfied 
if the quartic symmetroid is defined by a {\sl regular} web of quadrics in 
$\PP^3$, see \cite[Definition 2.1.2]{Cos}.
\hfill\par\noindent 2) Coble asserts in \cite{Coble1} that the 
converse to Proposition \ref{Coble} should also be true, presumably 
under suitable generality assumptions. This also relates to the 
question mentioned in \cite{EiPo} of describing when a collection 
of $10$ points in $\PP^3$ are determinantal.
\end{remark}

\section{Further results}
\label{How-many}

One may now ask which intersection numbers can actually occur and in which
geometric situation this can happen. We are far from having a complete 
answer to this question, but want to state a number of results in this direction. 

We start by considering Veronese surfaces which intersect in $10$ points. 
It is easy to find examples of surfaces $X_1$ and $X_2$ intersecting transversally 
in $10$ points. For this, one can start with an arbitrary surface 
$X_1\subset\operatorname{Gr}(\PP^1,\PP^3)\subset\PP^5$. For a general automorphism 
$\varphi$ of $\PP^3$ the surface $X_2=\varphi(X_1)$ intersects 
$X_1$ transversally by Kleiman's transversality theorem and since 
both $X_1$ and $X_2$ have the same cohomology class we have $X_1.X_2=10$.
Actually, we have the following

\begin{proposition}\label{prop3.2}
Let $X_1, X_2$ be two Veronese surfaces in $\PP^5$ 
intersecting in $10$ points. Then $X_1\cup X_2$ is contained in a 
hyperquadric $Q$ and one of the following cases
occurs:
\begin{enumerate}
\item[\rm{(i)}]
$Q$ has rank $6$ and $X_1$ and $X_2$ lie in the same cohomology class,
\item[\rm{(ii)}]
The rank of $Q$ is $4$ and $X_1\cap X_2\cap\operatorname{Sing} Q\neq\emptyset$,
\end{enumerate}
\end{proposition}

\begin{Proof}
The existence of $Q$ follows from Proposition \ref{prop1.4}. If $Q$ has rank
$6$ then $X_1$ and $X_2$ must have the same class by the proof of Proposition
\ref{prop2.1}. The case rank $Q=5$ is excluded by Proposition \ref{prop2.2}
and the case of rank $Q\le 4$ and $X_1\cap X_2\cap\operatorname{Sing}
Q=\emptyset$ is excluded by the remark at the end of Section 
\ref{Veronese-on-hyperquadrics}. We are now left to exclude the
case where rank $Q=3$ and $X_1\cap X_2\cap\operatorname{Sing} Q\neq
\emptyset$. We will make use of the diagram and the computations
at the end of the proof of Proposition \ref{prop2.4}.

For a Veronese surface $X\subset Q$ with rank $Q=3$ the class of $\tilde
X$ in $\tilde Q$ equals $2H^2$. The fibres of the map 
$\pi|_{\tilde X}:\tilde X\rightarrow C\cong \PP^1$ are conics and hence 
this linear system is a subsystem of $|2l|$ on $X\cong\PP^2$, 
i.e. contained in some system of the form $|2l-\sum\alpha_i P_i|$. 
We have $(2l-\sum\alpha_iE_i)^2=4-\sum\alpha^2_i=0$.
This implies either $\alpha_1=\ldots=\alpha_4=1$ or $\alpha_1=2$. 
But by the argument in the proof of Proposition \ref{prop2.4} we
already know that $X$ intersects the vertex of $Q$ in a finite 
non-empty set of points. Since the fibres of the map $\pi|_{\tilde X}$ 
are conics the first case can only occur if we have the linear system of
conics through $4$ points in general position. In this case 
$\tilde X$ is mapped to a $\PP^1$ and the general fibre is an irreducible conic, 
which contradicts what we have. This implies that the linear system is given 
by $|2l-2P|$, in particular $X$ meets the vertex of $Q$ in exactly one point $P$ 
and that this intersection is not-transversal.

Returning to the case we have to exclude, 
we can assume that $X_1$ is given by the $2\times 2$-minors of the matrix
$$
\begin{pmatrix}
x_0 & x_1 & x_2\\
x_1 & x_3 & x_4\\
x_3 & x_4 & x_5
\end{pmatrix}
$$
and so a typical rank $3$ hyperquadric through $X_1$ is $Q=\{x_0x_3-x^2_1=0\}$.
The vertex $V$ of $Q$ is the plane $\{x_0=x_1=x_3=0\}$. The intersection of $X_1$ and
$V$ is then defined by the ideal $(x_2x_4,x^2_2,x^2_4)$, i.e. the
first infinitesimal neighborhood of $P$. The same holds for the second 
Veronese surface $X_2\subset Q$. Since both surfaces $\tilde X_1$ and
$\tilde X_2$ have class $2H^2$ on $\tilde Q$, it follows that
$\tilde X_1$ and $\tilde X_2$ must meet in $6$ other points (different of $P$).
If these points are all in different fibers of $\pi$, then
 $X_1$ and $X_2$ meet in a (non-reduced) scheme of length 
$9(=6+3)$ and thus we got a contradiction. Otherwise, since
all the fibers of $\pi|_{\tilde X_i}$ are conics (whose images in $\PP^5$
already go through the fixed point $P$), the $2$-planes spanned by the 
fibres of $\pi|_{\tilde X_1}$  and $\pi|_{\tilde X_2}$ over some point
in $C$ must coincide. But then, by Bezout, the two conic fibres must intersect
in $4$ points, from which we conclude $X_1.X_2\ge 11$, again impossible.
This concludes the proof of the proposition.
\hfill\end{Proof}

\begin{remark}
As corollary of the proof of the previous proposition, we find out that
two Veronese surfaces in $\PP^5$ may intersect in a non-reduced scheme 
of length $9$. This is indeed the case, for two (general) 
Veronese surfaces $X_1$ and $X_2$ lying on a rank $3$ hyperquadric $Q$ such that 
$X_1\cap X_2\cap\operatorname{Sing}Q$ consists of a single point $P$. 
As above, both surfaces cut out the first infinitesimal neighborhood of $P$ 
on the vertex of the hyperquadric $Q$ and meet further in $6$ simple points.
In particular their intersection is non-reduced.
\end{remark}

However, we show that

\begin{proposition}
Two Veronese surfaces in $\PP^5$ cannot intersect transversally in $9$
points.
\end{proposition}

\begin{Proof}
Assume that there exist Veronese surfaces $X_1$ and $X_2$ in $\PP^5$
such that $W=X_1\cap X_2$ is a reduced set of length $9$. We consider
again $W$ as a subscheme of $X_1\cong\PP^2$ and we recall from Lemma
\ref{five-reg} that $h^1({\cal I}_W(4))=0$.  Equivalently, this means
that the linear system $\delta:=|4l-W|$ on $X_1\cong\PP^2$ has
projective dimension $5$. It will then be enough to show
that the linear system $|4l-W|$ of plane quartics through $W$ contains
a smooth curve $C$, since in this case the restriction of $|4l-W|$
defines, via taking the residual intersection, a $g^{4+\varepsilon}_7$
on $C$ with $\varepsilon\ge 0$, thus contradicting Clifford's theorem.

In order to show the existence of such a $C$ we consider the surface
$S$ given by blowing up the set $W$ on $X_1$. It will be enough to
check that the linear system $|4l-W|$ is base point free and defines a
morphism $S\rightarrow S'\subset\PP^5$ whose image $S'$ has no worse
than isolated singularities. We want to do this using Reider's theorem
(see \cite[Theorem 2.1]{BS} in characteristic $0$, and \cite{SB},
\cite{Nak}, \cite[Theorem 2.4]{Ter} in positive characteristic) and
for this purpose we write $4 l-W=L+K_{\delta}$ where
$K_{\delta}=-3l+W$ is the canonical divisor on $S$ while $L=7l-2W$. In
order to apply a Reider type theorem we need to check that $L^2\ge 9$
and that $L$ is nef and big. The first is clear since
$L^2=49-36=13$. We do not know that $L$ is nef and big, but in the
proof of Reider's theorem (as in the proof of its positive
characteristic counterparts, cf \cite{Ter}) this assumption is only
used to conclude that $h^1(K_{\delta}+L)=0$, which we already know
since $h^1(\calI_W(4))=0$ by Lemma \ref{five-reg}.

We may now argue as follows. If $|4l-W|=|K_{\delta}+L|$ is not
base-point free, respectively very ample, then there exists a curve a
$D$ such that $L-2D$ is $\QQ$-effective and such that
$$
D^2\ge L.D-k-1
$$ where $k=0$, respectively $1$. We write
$D=al-\sum\limits^9_{i=1}b_i E_i$. Since $L-2D$ must be
$\QQ$-effective we see immediately that $a\le 3$. For $a=3$ we obtain
$$
g-\sum\limits^9_{i=1}b^2_i \ge 21-2\sum\limits^9_{i=1} b_i-k-1
$$
respectively
$$
\sum\limits^9_{i=1}b_i(b_i-1)\le-2+k
$$
which gives a contradiction. For $a=1,2$ the same calculation gives
$$
\sum\limits^9_{i=1}(b_i-1)^2\le 4+k\qquad (a=1),
$$
respectively
$$
\sum\limits^9_{i=1}(b_i-1)^2\le k.
$$ On the other hand, since the quadrics through $X_2$ cut out $X_2$,
then at most $4$ points of $W\subset\PP^2$ can be collinear and at
most $8$ points of $W$ can lie on the same conic. This shows that
$|4l-W|$ is base point free on $S$ and that $S'$ has at most isolated
singularities, which is our claim, and this concludes the proof of the
proposition.
\end{Proof}

\begin{remark}
One may construct pairs of Veronese surfaces on suitable 
smooth or nodal cubic hypersurfaces in $\PP^5$ which meet 
in $1$, $2$, $3$, $5$ or $6$ simple points. It is also
possible to check in {\sl Macaulay} \cite{Macaulay} that if 
$X\subset\PP^5$ is a Veronese surface and $\varphi$ is a general 
linear automorphism of $\PP^5$ fixing 
$m\in\{1,2,3,5\}$ (general) points on $X$, then $X$
and $\varphi(X)$ meet exactly at those $m$ points.
\end{remark}

\bibliographystyle{alpha}

\bigskip
\noindent
\begin{tabular}{lll}
D. Eisenbud & K. Hulek & S. Popescu\\
Dept. of Mathematics & Inst. f\"ur Mathematik & Dept. of Mathematics\\
UC Berkeley & Univ. Hannover & SUNY Stony Brook\\
Berkeley CA 94720 & D 30060 Hannover & Stony Brook, NY 11794-3651\\
USA & Germany & USA\\
\\
{\tt de@msri.org} & {\tt hulek@math.uni-hannover.de} & {\tt sorin@math.sunysb.edu}\\
\end{tabular}
\end{document}